\nonstopmode
\documentclass[a4paper, 10pt]{amsart}
\usepackage{latexsym}
\usepackage{fancyhdr}
\usepackage{amsmath, amssymb}
\usepackage[ansinew]{inputenc}

\theoremstyle{plain}
\newtheorem*{definition*}{Definition}

\newtheorem*{lemma*}{Lemma}
\newtheorem{lemma}[subsection]{Lemma}
\newtheorem*{theorem*}{Theorem}
\newtheorem{theorem}[subsection]{Theorem}
\newtheorem*{proposition*}{Proposition}
\newtheorem{proposition}[subsection]{Proposition}
\newtheorem*{corollary*}{Corollary}

\theoremstyle{remark}
\newtheorem*{remark*}{Remark}
\newtheorem{remark}[subsection]{Remark}

\def\reg{\mathrm{reg}}

\title[Lifting smooth curves over invariants, II]
{Lifting smooth curves over invariants\\ 
for representations of compact Lie groups, II}

\author[A. Kriegl, M. Losik, P.W. Michor, A. Rainer]
{Andreas Kriegl, Mark Losik, Peter W. Michor, and Armin Rainer}

\begin{document}

\begin{abstract}
Any sufficiently often differentiable curve 
in the orbit space of a compact Lie group representation 
can be lifted to a once differentiable curve 
into the representation space. 
\end{abstract}

\thanks{M.L., P.W.M. and A.R. were supported by 
`Fonds zur F\"orderung der wissenschaftlichen Forschung, Projekt P 14195 MAT'}
\date {\today}

\maketitle

\section{Introduction}

In \cite{rep-lift} the following problem was investigated. 
Consider an orthogonal representation of a compact Lie group $G$ 
on a real finite dimensional Euclidean vector space $V$. 
Let $\sigma_1,\ldots,\sigma_n$ be a system of homogeneous generators 
for the algebra $\mathbb{R}[V]^G$ of invariant polynomials on $V$. 
Then the mapping 
$\sigma = (\sigma_1,\ldots,\sigma_n) : V \rightarrow \mathbb{R}^n$ 
induces a homeomorphism between the orbit space $V/G$ 
and the semialgebraic set $\sigma(V)$. 
Suppose a smooth curve 
$c : \mathbb{R} \rightarrow V/G = \sigma(V) \subseteq \mathbb{R}^n$ 
in the orbit space is given (smooth as curve in $\mathbb{R}^n$), 
does there exist a smooth lift to $V$, 
i.e., a smooth curve $\bar{c} : \mathbb{R} \rightarrow V$ with 
$c = \sigma \circ \bar{c}$\,?

It was shown in \cite{rep-lift} that a real analytic curve in $V/G$ 
admits a local real analytic lift to $V$, and that  
a smooth curve in $V/G$ admits a global smooth lift, 
if certain genericity conditions are satisfied.
In both cases the lifts 
may be chosen orthogonal to each orbit they meet and then they are  
unique up to a transformation in $G$, 
whenever the representation of $G$ on $V$ is polar, 
i.e., admits sections. 

In this paper we treat the same problem under weaker differentiability 
conditions for $c : \mathbb{R} \rightarrow V/G$ and without 
the mentioned genericity conditions. 
In section 3 we show that a continuous curve in the orbit space $V/G$ 
allows a global continuous lift to $V$. 
As a consequence we can prove in section 4 that a sufficiently often 
differentiable curve in $V/G$ can be lifted to a once differentiable 
curve in $V$. 
What we mean by sufficiently often differentiable will be 
specified there.
 
In the special case that the symmetric group $S_n$ 
is acting on $\mathbb{R}^n$, 
in other words (see \cite{rep-lift}), if 
smooth parameterizations of the roots of smooth curves of polynomials 
with all roots real are looked for, 
the following results were proved in \cite{rootsII}: 
Any differentiable lift of a $C^{2 n}$-curve (of polynomials) 
$c : \mathbb{R} \rightarrow \mathbb{R}^n/S_n$ is actually $C^1$, 
and there always exists a twice differentiable 
but in general not better lift of $c$, 
if it is of class $C^{3 n}$. 
Note that here the differentiability assumptions on $c$ are not the
weakest possible which is shown by the case $n=2$, 
elaborated in \cite{roots} 2.1. 
The proof there is based on the fact that the roots of a 
$C^n$-curve of 
polynomials $c : \mathbb{R} \rightarrow \mathbb{R}^n/S_n$ 
may be chosen differentiable with locally bounded derivative; 
this is due to Bronshtein \cite{bronshtein} and 
Wakabayashi \cite{wakabayashi}. 
Therefore, our long-term objective is to prove the existence of a 
twice differentiable lift also in the general setting. 
The key is the generalization of Bronshtein's and Wakabayashi's result 
which seems to be difficult.  

The polynomial results have applications in the theory 
of partial differential equations and perturbation 
theory, see \cite{perturb}.

\section{Preliminaries}

\subsection{The setting}

Let $G$ be a compact Lie group and let $\rho : G \rightarrow O(V)$ be an 
orthogonal representation in a real finite dimensional Euclidean 
vector space $V$ with inner product $\langle \quad |\quad \rangle$. 
By a classical theorem of Hilbert and Nagata, 
the algebra $\mathbb{R}[V]^{G}$ of invariant polynomials on $V$ 
is finitely generated. 
So let $\sigma_1,\ldots,\sigma_n$ be a system of homogeneous generators 
of $\mathbb{R}[V]^{G}$ of positive degrees $d_1,\ldots,d_n$. 
We may assume that $\sigma_1 : v \mapsto \langle v|v \rangle$ 
is the inner product. 
Consider the \emph{orbit map} 
$\sigma = (\sigma_1,\ldots,\sigma_n) : V \rightarrow \mathbb{R}^n$. 
Note that, 
if $(y_1,\ldots,y_n) = \sigma(v)$ for $v \in V$, then 
$(t^{d_1} y_1,\ldots,t^{d_n} y_n) = \sigma(t v)$ for $t \in \mathbb{R}$, 
and that $\sigma^{-1}(0) = \{0\}$. 
The image $\sigma(V)$ is a semialgebraic set in the categorical quotient 
$V/\!\!/G := \{y \in \mathbb{R}^n : P(y) = 0 ~\mbox{for all}~ P \in I\}$ 
where $I$ is the ideal of relations between $\sigma_1,\ldots,\sigma_n$. 
Since $G$ is compact, $\sigma$ is proper and separates orbits of $G$, 
it thus induces a homeomorphism between $V/G$ and $\sigma(V)$.

\subsection{The slice theorem} \label{slice}
For a point $v \in V$ we denote by $G_v$ its isotropy group and by 
$N_v = T_v(G.v)^{\bot}$ the normal subspace of the orbit $G.v$ at $v$. 
It is well known that there exists a $G$-invariant open neighborhood $U$ 
of $v$ which is real analytically $G$-isomorphic to the crossed product 
(or associated bundle) $G \times_{G_v} S_v = (G \times S_v)/G_v$, 
where $S_v$ is a ball in $N_v$ with center at the origin. 
The quotient $U/G$ is homeomorphic to $S_v/G_v$. 
It follows that the problem of local lifting curves in $V/G$ passing 
through $\sigma(v)$ reduces to the same problem for curves in $N_v/G_v$ 
passing through $0$. 
For more details see \cite{rep-lift}, \cite{luna} and \cite{schwarz1}, 
theorem 1.1.

A point $v \in V$ (and its orbit $G.v$ in $V/G$) is called regular 
if the isotropy  representation  
$G_v \rightarrow O(N_v)$  is trivial.
Hence a neighborhood of this point is analytically $G$-isomorphic to 
$G/G_v \times S_v \cong G.v \times S_v$. 
The set $V_{\reg}$ of regular points is open and dense in $V$, 
and the projection $V_{\reg} \rightarrow V_{\reg}/G$ 
is a locally trivial fiber bundle. 
A non regular orbit or point is called singular. 

\subsection{Removing fixed points} \label{fix}

Let $V^G$ be the space of $G$-invariant vectors in $V$, 
and let $V'$ be its orthogonal complement in $V$. 
Then we have $V = V^G \oplus V'$, 
$\mathbb{R}[V]^G = \mathbb{R}[V^G] \otimes \mathbb{R}[V']^G$ and 
$V/G = V^G \times V'/G$. 

\begin{lemma*}
Any lift $\bar{c}$ of a curve $c=(c_0,c_1)$ of class $C^k$ 
($k = 0,1,\ldots,\infty,\omega$) in $V^G \times V'/G$ has the form 
$\bar{c} = (c_0,\bar{c}_1)$, 
where $\bar{c}_1$ is a lift of $c_1$ to $V'$ of class 
$C^k$ ($k = 0,1,\ldots,\infty,\omega$). 
The lift $\bar{c}$ is orthogonal if and only if the lift 
$\bar{c}_1$ is orthogonal. \qed
\end{lemma*}

\subsection{Multiplicity} \label{mult}

For a continuous function $f$ defined near $0$ in $\mathbb{R}$, 
let the \emph{multiplicity} or \emph{order of flatness} $m(f)$ at 
$0$ be the supremum of all integers $p$ such that $f(t) = t^p g(t)$ near 
$0$ for a continuous function $g$. 
If $f$ is $C^n$ and $m(f) < n$, 
then $f(t) = t^{m(f)} g(t)$, 
where now $g$ is $C^{n-m(f)}$ and $g(0) \ne 0$. 
Similarly, 
one can define multiplicity of a function at any $t \in \mathbb{R}$.

\begin{lemma*}
Let $c = (c_1,\ldots,c_n)$ be a curve in $\sigma(V) \subseteq \mathbb{R}^n$, 
where $c_i$ is $C^{d_i}$, 
for $1 \le i \le n$, 
and $c(0) = 0$. 
Then the following two conditions are equivalent:
\begin{enumerate}
\item $c_1(t) = t^2 c_{1,1}(t)$ near $0$ for a continuous function $c_{1,1}$;
\item $c_i(t) = t^{d_i} c_{i,i}(t)$ near $0$ for a continuous function $c_{i,i}$, 
for all $1 \le i \le n$.
\end{enumerate}
\end{lemma*}

\proof
The proof of the nontrivial implication $(1) \Rightarrow (2)$ 
is the same as in the smooth case with $r = 1$, 
see \cite{rep-lift} 3.3. for details. \qed

\section{Lifting continuous curves over invariants}

\begin{proposition}      \label{cont}
Let 
$c = (c_1,\ldots,c_n) : \mathbb{R} \rightarrow V/G 
= \sigma(V) \subseteq \mathbb{R}^n$ 
be continuous. 
Then there exists a global continuous lift 
$\bar{c} : \mathbb{R} \rightarrow V$ of $c$. 
\end{proposition}

This result is due to Montgomery and Yang \cite{mont} 
see also \cite{bredon}. 
We present a short proof adapted to our setting:

\proof  

We will make induction on the size of $G$. 
More precisely, 
for two compact Lie groups $G'$ and $G$ we denote $G' < G$, 
if 
\begin{itemize}
\item $\dim G' < \dim G$ \qquad or
\item if $\dim G' = \dim G$, 
      then $G'$ has less connected components than $G$ has.
\end{itemize}
In the simplest case, 
when $G = \{e\}$ is trivial, 
we find $\sigma(V) = V/G = V$, 
whence we can put $\bar{c} := c$.

Let us assume that for any $G' < G$ and any continuous 
$c : \mathbb{R} \rightarrow V/G'$ there exists a global continuous lift 
$\bar{c} : \mathbb{R} \rightarrow V$ of $c$, 
where $G' \rightarrow O(V)$ is an orthogonal representation on 
an arbitrary real finite dimensional Euclidean vector space $V$.

We shall prove that then the same is true for $G$. 
Let $c : \mathbb{R} \rightarrow V/G = \sigma(V) \subseteq \mathbb{R}^n$ 
be continuous.  By lemma \ref{fix}, 
we may remove the nontrivial fixed points of the $G$-action on $V$ and 
suppose that $V^{G} = \{0\}$. 
The set $c^{-1}(0)$ is closed in $\mathbb{R}$ and,  consequently, 
$c^{-1}(\sigma(V) \backslash \{0\}) = \mathbb{R} \backslash c^{-1}(0)$ 
is open in $\mathbb{R}$. 
Thus, we can write 
$c^{-1}(\sigma(V) \backslash \{0\}) = \bigcup_{i \in I} (a_i,b_i)$, 
a disjoint union, where $a_i,b_i \in \mathbb{R} \cup \{\pm \infty\}$ 
with $a_i < b_i$ such that each $(a_i,b_i)$ is maximal with respect 
to not containing zeros of $c$, 
and $I$ is an at most countable set of indices. 
In particular, 
we have $c(a_i) = c(b_i) = 0$ for all $a_i, b_i \in \mathbb{R}$ 
appearing in the above presentation.

We assert that on each $(a_i,b_i)$ there exists a continuous lift 
$\bar{c} : (a_i,b_i) \rightarrow V \backslash \{0\}$ 
of the restriction 
$c|_{(a_i,b_i)} : (a_i,b_i) \rightarrow \sigma(V) \backslash \{0\}$. 
In fact, 
since $V^G = \{0\}$, 
for all $v \in V \backslash \{0\}$ the isotropy groups $G_v$, 
acting orthogonally on $N_v$, 
satisfy $G_v < G$. 
Therefore, 
by induction hypothesis and by \ref{slice}, 
we find local continuous lifts of $c|_{(a_i,b_i)}$ near any 
$t \in (a_i,b_i)$ and through all $v \in \sigma^{-1}(c(t))$. 
Suppose $\bar{c}_1 : (a_i,b_i) \supseteq (a,b) \rightarrow V \backslash \{0\}$ 
is a local continuous lift of $c|_{(a_i,b_i)}$ with maximal domain $(a,b)$, 
where, say, $b < b_i$. 
Then there exists a local continuous lift $\bar{c}_2$ 
of $c|_{(a_i,b_i)}$ near $b$, 
and there is a $t_0 < b$ such that both $\bar{c}_1$ and $\bar{c}_2$ 
are defined near $t_0$. 
Since $\bar{c}_1(t_0)$ and $\bar{c}_2(t_0)$ lie in the same orbit,  
there must exist a $g \in G$ such that 
$\bar{c}_1(t_0) = g . \bar{c}_2(t_0)$. 
But then, 
\[\bar{c}_{12}(t) := \left\{ \begin{array}{r@{\quad \mbox{for} \quad}l} 
  \bar{c}_1(t) & t \le t_0 \\ g . \bar{c}_2(t) & t \ge t_0  \end{array} \right.
\]
is a local continuous lift of $c|_{(a_i,b_i)}$ defined on a larger 
interval than $\bar{c}_1$. 
Thus we have shown that each local continuous lift of $c|_{(a_i,b_i)}$ 
defined on an open interval $(a,b) \subseteq (a_i,b_i)$ can be extended 
to a larger interval whenever $(a,b) \subsetneq (a_i,b_i)$. 
This proves the assertion.

We put $\bar{c}|_{c^{-1}(0)} := 0$, since, 
by $\sigma^{-1}(0) = \{0\}$, this is the only choice. 
Then $\bar c$ is also continuous at points $t_0\in c^{-1}(0)$ since
$\langle \bar{c}(t) | \bar{c}(t) \rangle = \sigma_1(\bar{c}(t)) = c_1(t)$ 
converges to $0$ as $t \to t_0$.
\qed

\section{Lifting differentiably}

Throughout the whole section we let $d\ge 2$
be the maximum of all degrees of systems of minimal generators of 
invariant polynomials of all slice representations of $\rho$. 
Of these there are only finitely many isomorphism types. 

\begin{lemma}           \label{reg}
A curve
$c : \mathbb{R} \rightarrow V/G = \sigma(V) \subseteq \mathbb{R}^n$
of class $C^d$ admits an orthogonal $C^d$-lift $\bar{c}$ in a
neighborhood of a regular point $c(t_0) \in V_{\reg}/G$.
It is unique up to a transformation from $G$.
\end{lemma}

\proof

The proof works analogously as in the smooth case,
see \cite{rep-lift} 3.1. \qed

\begin{theorem}  \label{diff}
Let
$c = (c_1,\ldots,c_n) : \mathbb{R} \rightarrow V/G
= \sigma(V) \subseteq \mathbb{R}^n$ be a curve of class $C^d$.
Then for any $t_0 \in \mathbb{R}$ there exists a local lift
$\bar{c}$ of $c$ near $t_0$ which is differentiable at $t_0$.
\end{theorem}

\proof

We follow partially the algorithm given in \cite{rep-lift} 3.4.
Without loss of generality we may assume that $t_0 = 0$.
We show the existence of local lifts of $c$ which are differentiable
at $0$ through any $v \in \sigma^{-1}(c(0))$.
By lemma \ref{fix}
we can assume $V^G =\{0\}$.

If $c(0) \ne 0$ corresponds to a regular orbit,
then unique orthogonal $C^d$-lifts defined near $0$ exist through
all $v \in \sigma^{-1}(c(0))$,
by lemma \ref{reg}.

If $c(0) = 0$,
then $c_1$ must vanish of at least second order at $0$,
since $c_1(t) \ge 0$ for all $t \in \mathbb{R}$.
That means $c_1(t) = t^2 c_{1,1}(t)$ near $0$ for a continuous function
$c_{1,1}$ since $c_1$ is $C^2$.
By the multiplicity lemma \ref{mult}
we find that $c_i(t) = t^{d_i} c_{i,i}(t)$ near $0$ for $1 \le i \le n$,
where $c_{1,1},c_{2,2},\ldots,c_{n,n}$ are continuous functions.
We consider the following curve in $\sigma(V)$ which is continuous since
$\sigma(V)$ is closed in $\mathbb{R}^n$,
see \cite{procesischwarz}:
\begin{eqnarray*}
c_{(1)}(t) :&=& (c_{1,1}(t),
c_{2,2}(t),\ldots,c_{n,n}(t))\\ &=&
(t^{-2} c_1(t),t^{-d_2}c_2(t),\ldots,t^{-d_n} c_n(t)).
\end{eqnarray*}
By proposition \ref{cont},
there exists a continuous lift $\bar{c}_{(1)}$ of $c_{(1)}$.
Thus, $\bar{c}(t) := t \cdot \bar{c}_{(1)}(t)$ is a local lift of $c$
near $0$ which is differentiable at $0$:
\[\sigma(\bar{c}(t)) = \sigma(t \cdot \bar{c}_{(1)}(t))
= (t^2 c_{1,1}(t),\ldots,t^{d_n} c_{n,n}(t)) = c(t),\]
and
\[\lim_{t \rightarrow 0} \frac{t \cdot \bar{c}_{(1)}(t)}{t}
=  \lim_{t \rightarrow 0} \bar{c}_{(1)}(t) = \bar{c}_{(1)}(0).\]
Note that $\sigma^{-1}(0) = \{0\}$,
therefore we are done in this case.

If $c(0) \ne 0$ corresponds to a singular orbit,
let $v$ be in $\sigma^{-1}(c(0))$ and consider the isotropy representation
$G_v \rightarrow O(N_v)$.
By \ref{slice},
the lifting problem reduces to the same problem for curves in $N_v/G_v$
now passing through $0$.
\qed

\begin{lemma}    \label{topolemma}
Consider a continuous curve $c : (a,b) \to X$ in a compact metric space $X$.
Then the set $A$ of all accumulation points of $c(t)$ as
$t \searrow a$ is connected.
\end{lemma}

\proof
On the contrary suppose that $A = A_1 \cup A_2$,
where $A_1$ and $A_2$ are disjoint open and closed subsets of $A$.
Since $A$ is closed in $X$,
also $A_1$ and $A_2$ are closed in $X$.
There exist disjoint open subsets $A_1',A_2' \subseteq X$ with
$A_1 \subseteq A_1'$ and $A_2 \subseteq A_2'$.
Consider $F := X \backslash (A_1' \cup A_2')$ which is closed in $X$
and hence compact.
Since $c$ visits $A_1'$ and $A_2'$ infinitely often and
$c^{-1}(A_1')$ and $c^{-1}(A_2')$ are disjoint and open in $\mathbb R$, 
there exists a
sequence $t_m \to a$ and $c(t_m) \in F$
for all $m$.
By compactness of $F$,
this sequence has a cluster point $y$ in $F$.
Hence $y$ is in $A$ by definition, which contradicts 
$F \cap A = \emptyset$. \qed

\begin{theorem}  \label{gldiff}
Let
$c = (c_1,\ldots,c_n) : \mathbb{R} \rightarrow V/G
  = \sigma(V) \subseteq \mathbb{R}^n$
be a curve of class $C^d$.
Then there exists a global differentiable lift
$\bar{c} : \mathbb{R} \rightarrow V$ of $c$.
\end{theorem}

\proof

The proof,
as the one of proposition \ref{cont},
will be carried out by induction on the size of $G$.

If $G = \{e\}$ is trivial,
then $\bar{c} := c$ is a global differentiable lift.

So let us assume that for any $G' < G$ and any
$c : \mathbb{R} \rightarrow V/G'$
satisfying the differentiability conditions of the theorem
there exists a global differentiable lift
$\bar{c} : \mathbb{R} \rightarrow V$ of $c$,
where $G' \rightarrow O(V)$ is an orthogonal representation on an
arbitrary real finite dimensional Euclidean vector space $V$.

We shall prove that the same is true for $G$.
Let
$c = (c_1,\ldots,c_n) : \mathbb{R} \rightarrow V/G
  = \sigma(V) \subseteq \mathbb{R}^n$ be of class $C^d$.
We may assume that $V^G = \{0\}$,
by lemma \ref{fix}.
As in the proof of proposition \ref{cont} we can write
$c^{-1}(\sigma(V) \backslash \{0\}) = \bigcup_{i} (a_i,b_i)$,
a disjoint union,
where $a_i,
b_i \in \mathbb{R} \cup \{\pm \infty\}$ with
$a_i < b_i$.
In particular,
we have $c(a_i) = c(b_i) = 0$ for all $a_i,
b_i \in \mathbb{R}$ appearing in the above presentation.

{\it Claim: On each $(a_i,b_i)$ there exists a differentiable lift
$\bar{c} : (a_i,b_i) \rightarrow V \backslash \{0\}$ of the restriction
$c|_{(a_i,b_i)} : (a_i,b_i) \rightarrow \sigma(V) \backslash \{0\}$}.
The lack of nontrivial fixed points guarantees that for all
$v \in V \backslash \{0\}$ the isotropy groups $G_v$ acting on
$N_v$ satisfy $G_v < G$.
Therefore, by induction hypothesis and by \ref{slice},
we find local differentiable lifts of $c|_{(a_i,b_i)}$ near any
$t \in (a_i,b_i)$ and through all $v \in \sigma^{-1}(c(t))$.
Suppose that
$\bar{c}_1 : (a_i,b_i) \supseteq (a,b) \rightarrow V \backslash \{0\}$
is a local differentiable lift of $c|_{(a_i,b_i)}$ with maximal domain $(a,b)$,
where,
say, $b < b_i$.
Then there exists a local differentiable lift $\bar{c}_2$ of
$c|_{(a_i,b_i)}$ near $b$,
and there exists a $t_0 < b$ such that both $\bar{c}_1$ and
$\bar{c}_2$ are defined near $t_0$.
We may assume without loss that $\bar{c}_1(t_0) = \bar{c}_2(t_0) =: v_0$,
by applying a transformation $g \in G$ to $\bar{c}_2$,
say.
We want to show that we can arrange the lift $\bar{c}_2$
in such a way that its derivative at $t_0$ matches with the derivative
of $\bar{c}_1$ at $t_0$.
We decompose $\bar c_i'(t_0)=\bar c_i'(t_0)^\top + \bar c_i'(t_0)^\bot$ into the parts 
tangent to the orbit $G.v_0$ and normal to it. 

First we deal with the normal parts $\bar c_i'(t_0)^\bot\in V$. 
We consider the projection
$p : G.S_{v_0} \cong G \times_{G_{v_0}} S_{v_0} \to G/G_{v_0} \cong G.v_0$
of the fiber bundle associated to the principal bundle
$\pi : G \to G/G_{v_0}$.
Then, for $t$ close to $t_0$, $\bar{c}_1$ and $\bar{c}_2$
are differentiable curves in $G.S_{v_0}$,
whence $p \circ \bar{c}_i$ $(i =1,2)$ are differentiable curves in
$G/G_{v_0}$ which admit differentiable lifts $g_i$ into $G$ with
$g_i(t_0) = e$ (via the horizontal lift of a principal connection,
say).
Consequently, $t \mapsto g_i(t)^{-1}.\bar{c}_i(t)=:\tilde c_i(t)$ 
are differentiable
lifts of $c|_{(a_i,b_i)}$ near $t_0$ which lie in $S_{v_0}$,
whence
$\tilde c_i'(t_0)=\left.\frac{d}{d t}\right|_{t = t_0} (g_i(t)^{-1}.\bar{c}_i(t))
  = - g_i'(t_0).v_0 + \bar{c}_i'(t_0) \in N_{v_0}$.
So, $\bar{c}_i'(t_0)^\top = (g_i'(t_0).v_0)^\top = g_i'(t_0).v_0$,
and so for the normal part  we get 
$\bar{c}_i'(t_0)^\bot =\tilde c_i'(t_0)$.

Since $\tilde c_i$ lie in $S_{v_0}$ 
we can change to the isotropy representation
$G_{v_0} \rightarrow O(N_{v_0})$ (using the same letters $\sigma_i$ 
for the generators of $\mathbb R[N_{v_0}]^{G_{v_0}}$).
We can suppose that $v_0 = 0$, i.e., $c(t_0) = 0$.

Recall the continuous curve in $\sigma(V)$
defined in the proof of theorem \ref{diff} which depends on the point $t_0$:
\[c_{(1,t_0)}(t)
  := ((t-t_0)^{-2} c_1(t), (t-t_0)^{-d_2} c_2(t),\ldots,(t-t_0)^{-d_n} c_n(t)).\]
We find that for $i = 1,2$:
\[\sigma(\tilde{c}_i'(t_0))
= \sigma \left(\lim_{t \rightarrow t_0} \frac{\tilde{c}_i(t)
  - \tilde{c}_i(t_0)}{t-t_0}\right)
= \lim_{t \rightarrow t_0} \sigma \left(\frac{\tilde{c}_i(t)}{t-t_0}\right)
= c_{(1,t_0)}(t_0).\]
So $\tilde{c}_1'(t_0)$ and $\tilde{c}_2'(t_0)$ are lying in the same orbit.
This shows also that
\begin{list}{$\blacklozenge$}{}
\item {} for any two lifts of $c$ near $t_0\in c^{-1}(0)$ 
      which are one-sided differentiable at $t_0$ the derivatives at 
      $t_0$ lie in the same $G$-orbit.
\end{list}
Thus,
there must exist a $g_0 \in G_{v_0}$ such that
$\bar c_1'(t_0)^\bot =\tilde{c}_1'(t_0) 
= g_0.\tilde{c}_2'(t_0)=g_0.\bar c_2'(t_0)^\bot=(g_0.\bar c_2)'(t_0)^\bot$.

Now we deal with the tangential parts.
We search for a differentiable curve $t \mapsto g(t)$ in $G$
with $g(t_0) = g_0$ and
\[\bar{c}_1'(t_0)^\top
  = \bigl(\tfrac{d}{d t}|_{t =t_0} (g(t).\bar{c}_2(t)) \bigr)^\top
  = g'(t_0).v_0 + g_0.\bar{c}_2'(t_0)^\top.\]
But this linear equation can be solved for $g'(t_0)$,
and, hence,
the required curve $t \mapsto g(t)$ exists.
Note that the normal parts still fit since
\[
\bigl(\tfrac{d}{d t}|_{t =t_0} (g(t).\bar{c}_2(t)) \bigr)^\bot
=\bigl(g'(t_0).v_0+g_0.\bar{c}_2'(t_0) \bigr)^\bot
=0+g_0.\bar{c}_2'(t_0)^\bot =\bar c_1'(t_0)^\bot.
\]
The two lifts $\bar c_1$ for $t\le t_0$ and $g.\bar c_2$ for $t\ge t_0$
fit together differentiably at $t_0$.
This proves the claim.

Now let $\bar{c} : (a_i,b_i) \rightarrow V \backslash \{0\}$
be the differentiable lift of $c|_{(a_i,b_i)}$ constructed above.
For $a_i \ne - \infty$, we put $\bar{c}(a_i) := 0$,
the only choice.
Consider the expression $\gamma(t) := \frac{\bar{c}(t)}{t-a_i}$
which is a differentiable curve in $V \backslash \{0\}$ for $t \in (a_i,b_i)$.
We want to show that $\lim_{t \searrow a_i} \gamma(t)$ exists.
For $t$ sufficiently close to $a_i$ we have
\[\sigma(\gamma(t)) = \sigma \left(\frac{\bar{c}(t)}{t-a_i}\right)
= c_{(1,a_i)}(t) \to c_{(1,a_i)}(a_i) \qquad \mbox{as} ~ t \searrow a_i,\]
where now
$c_{(1,a_i)}(t) := ((t-a_i)^{-2} c_1(t),
  (t-a_i)^{-d_2} c_2(t),\ldots,(t-a_i)^{-d_n} c_n(t))$.
Let $\bar{c}_{(1,a_i)}$ be a corresponding continuous lift of
$c_{(1,a_i)}$ which exists by proposition \ref{cont}.
This shows that the set $A$ of all accumulation points of
$\left( \gamma(t)\right)_{t \searrow a_i}$ lies in the orbit
$G.\bar{c}_{(1,a_i)}(a_i)$ through $\bar{c}_{(1,a_i)}(a_i)$.
By lemma \ref{topolemma}, $A$ is connected.
In particular,
the limit $\lim_{t \searrow a_i} \gamma(t)$ must exist,
if $G$ is a finite group.
In general let us consider the projection
$p : G.S_{v_1} \cong G \times_{G_{v_1}} S_{v_1} \to G/G_{v_1} \cong G.v_1$
of a fiber bundle associated to the principal bundle $\pi : G \to G/G_{v_1}$,
where we choose some $v_1 \in A$.
For $t$ close to $a_i$ the curve $t \mapsto \gamma(t)$ is differentiable
in $G.S_{v_1}$,
whence $t \mapsto p(\gamma(t))$ defines a differentiable curve in
$G/G_{v_1}$ which admits a differentiable lift $t \mapsto g(t)$ into $G$.
Now, $t \mapsto g(t)^{-1}.\gamma(t)$ is a differentiable curve in
$S_{v_1}$ whose accumulation points for $t \searrow a_i$ have to lie in
$G.v_1 \cap S_{v_1} = \{v_1\}$,
since $\sigma( g(t)^{-1}.\gamma(t)) = \sigma(\gamma(t))$.
That means that $t \mapsto g(t)^{-1}.\bar{c}(t)$ defines
a differentiable lift of $c|_{(a_i,b_i)}$,
for $t > a_i$ close to $a_i$,
whose one-sided derivative at $a_i$ exists:
\[\lim_{t \searrow a_i} \frac{g(t)^{-1}.\bar{c}(t)}{t-a_i}
= \lim_{t \searrow a_i} g(t)^{-1}.\gamma(t) = v_1.\]
Let $t\mapsto g(t)$ be extended smoothly to $(a_i,b_i)$ so that near $b_i$
it is constant  and replace $t\mapsto \bar c(t)$ by
$t\mapsto g(t)^{-1}\bar c(t)$.
Thus
\[\bar{c}'(a_i) := \lim_{t \searrow a_i} \frac{\bar{c}(t)}{t-a_i}
=v_1.\]
The same reasoning is true for $b_i \ne + \infty$.
Thus we have extended $\bar{c}$ differentiably to the closure of $(a_i,b_i)$.

Let us now construct a global differentiable lift of $c$ defined on the
whole of $\mathbb{R}$.
For isolated points $t_0 \in c^{-1}(0)$ the two differentiable lifts on
the neighboring intervals can be made to match differentiably, by applying
a fixed $g\in G$ to one of them by $\blacklozenge$.
Let $E$ be the set of accumulation points of $c^{-1}(0)$.
For connected components of $\mathbb R\setminus E$ we can proceed inductively
to obtain differentiable lifts on them.

We extend the lift by 0 on the set $E$ of
accumulation points of $c^{-1}(0)$.
Note that every
lift $\tilde c$ of $c$ has to vanish on $E$
and is continuous there since
$\langle \tilde c(t)|\tilde c(t) \rangle=\sigma_1(\tilde c(t))=c_1(t)$.
We also claim that any lift
$\tilde c$ of $c$ is differentiable at any point
$t'\in E$ with derivative 0.
Namely, the difference quotient $t\mapsto \frac{\tilde c(t)}{t-t'}$ is a lift
of the curve $c_{(1,t')}$
which vanishes at $t'$ by the following argument:
Consider the local lift $\bar{c}$ of $c$ near $t'$
which is differentiable at $t'$,
provided by theorem \ref{diff}.
Let $(t_m)_{m \in \mathbb{N}} \subseteq c^{-1}(0)$ be a sequence with
$t'\ne t_m \rightarrow t'$,
consisting exclusively of zeros of $c$.
Such a sequence always exists since $t'\in E$. 
Then we have
\[\bar{c}'(t')
= \lim_{t \rightarrow t'} \frac{\bar{c}(t) - \bar{c}(t')}{t-t'}
= \lim_{m \rightarrow \infty} \frac{\bar{c}(t_m)}{t_m -t'} = 0. \]
Thus $c_{(1,t')}(t')=\lim_{t\to t'}\sigma(\frac{\bar c(t)}{t-t'})
=\sigma(\bar c'(t'))=0$.
\qed

\remark \label{best}
Note that
the differentiability conditions of the curve $c$ in the current
section are best possible:
In the case when the symmetric group $S_n$ is acting in $\mathbb{R}^n$
by permuting the coordinates,
and $\sigma_1,\ldots,\sigma_n$ are the elementary
symmetric polynomials with degrees $1,\ldots,n$,
there need not exist a differentiable lift
if the differentiability assumptions made on $c$ are weakened,
see \cite{roots} 2.3. first example.
\endremark

\vspace{30pt}

\footnotesize

\noindent
{\scshape A. 
Kriegl: Institut f\"ur Mathematik, Universit\"at Wien, 
Nordbergstrasse~15, A-1090 Wien, Austria}

\emph{E-mail address:} 
{\ttfamily Andreas.Kriegl@univie.ac.at}

\vspace{10pt}

\noindent
{\scshape M. 
Losik: Saratov State University, 
ul. 
Astrakhanskaya, 83, 410026 Saratov, Russia}

\emph{E-mail address:} 
{\ttfamily losikMV@info.sgu.ru}

\vspace{10pt}

\noindent
{\scshape P.W. 
Michor: Institut f\"ur Mathematik, Universit\"at Wien, 
Nordbergstrasse~15, A-1090 Wien, Austria;} 
\emph{and:} 
{\scshape Erwin Schr\"odinger Institut f\"ur Mathematische Physik, 
Boltzmann\-gasse 9, A-1090 Wien, Austria}

\emph{E-mail address:} 
{\ttfamily Peter.Michor@esi.ac.at}

\vspace{10pt}

\noindent
{\scshape A. 
Rainer: Institut f\"ur Mathematik, Universit\"at Wien, 
Nordbergstrasse~15, A-1090 Wien, Austria}

\emph{E-mail address:} 
{\ttfamily armin\_rainer@gmx.net}

\end{document}